\documentclass[12pt,a4paper]{amsart}

\usepackage{amsaddr}
\usepackage{amssymb}
\usepackage[final]{showkeys}
\usepackage{microtype}
\usepackage{color}
\usepackage{graphicx}
\usepackage[hmargin=3cm,vmargin={3.5cm,4cm}]{geometry}

\numberwithin{equation}{section}

\allowdisplaybreaks

\begin{document}
	
	\title[]{On the role of Giovanni Giorgi in the history of operational methods for mathematical-physics problems}
	
	\author{Roberto Garra$^{1}$ \and Francesco Mainardi$^2$}
	\address{$^1$
		Section of Mathematics,  International Telematic University 
		\\Uninettuno,
		Via S. Martino della Battaglia  44, 00185 Roma, Italy.
		\\
		E-mail: roberto.garra@uninettunouniversity.net
		\\
		$^2$Department of Physics and Astronomy, University of Bologna, 
		\\ Via Irnerio 46, 40126 Bologna, Italy. 
		\\ E-mail: francesco.mainardi@unibo.it
		\\
	}

	\date{\today}
	
	\begin{abstract}
		In this paper we discuss the historical role played by Prof. Giovanni Giorgi (1871-1950) in the development of operational methods in mathematical physics. 
		In the literature, the analysis of the scientific contributions of Giorgi is discussed in many historical papers mainly about the $MKS\Omega$-system and its contributions in 
		the field of electrical engineering. Starting from the obituary written by Prof. Dario Graffi (1905-1990), here we analyze in detail the contributions given by Giorgi in mathematical physics, especially in the framework of the studies started by Heaviside (1850-1925) in order to obtain a symbolic representation of the solutions of differential equations emerging in physical problems. \\
		Moreover, we underline the little known contribution of Giorgi to the analysis of derivatives of any real order, namely fractional derivatives, according to the present notation. \\
		The main aim of this note is to underline the historical relevance of Giorgi in the mathematical foundations of operational methods in relation to the solutions of concrete problems emerging in applied physics.\\
		
		\smallskip
		
			\textbf{Keywords:} Giovanni Giorgi; operational calculus; history of fractional calculus
		\\	\textbf{MSC:} 01A70; 44A45
	\end{abstract}
	\maketitle
	
	\section{Introduction}  
	
	Giovanni Giorgi was born in Lucca on November 27, 1871.
	From an early age, he demonstrated a keen interest in scientific and technical matters. After completing his secondary education, he enrolled in 1888 at the University of Rome, where he studied under Luigi Cremona and Eugenio Beltrami and where he graduated as an engineer in 1893\footnote{This short introductory biography is mainly based on the Obituary by Prof. Graffi (1951). We underline that there is a wide literature about the life and work of Giovanni Giorgi. We refer for example to Rossi et al. (2002), D'Agostino and Rossi (2004), Frezza et al. (2015) and the references therein. }.
	He immediately devoted himself to his professional career, especially in the field of electrical engineering.
	\newpage
	From 1906 to 1921, he was director of the technological department of the Municipality of Rome. However, during this period, he also carried out significant scientific activity and taught several courses at the Faculty of Science and the School of Engineering of the University of Rome. 
	
	\smallskip
	
	In 1926, he was appointed to the chair of mathematical physics
	at the University of Cagliari. This fact is known to  surprise the community of physicists because the 
	committee preferred him as an engineer,  to Enrico Fermi, the young promising scientist who was enrolled as a teacher of mathematical physics at the University of Florence.
	It is known that in the same year Fermi got the chair of theoretical physics
	just at the University of Rome thanks to his mentor Prof. Orso  Corbino
	(Minister of Public Instruction) who was able to create for the first time in Italy the      chair of theoretical physics,
	see Rossi (2026).
	From this fact in Italy a dualism began between the chair of mathematical physics (administrated by mathematicians) 
	and the chair of theoretical physics (administrated by physicists): a long standing controversy that surely harmed the relationship between the departments of mathematics and physics all over Italy. In 1929 he moved to the same chair at the University of Palermo. From 1934, he was a teacher at the University of Rome as a full professor of electrical communications in the Faculty of Engineering. 
	
	\smallskip
	
	Having reached his age limits, he retired from official teaching in 1942 but continued his intense scientific and professional activity until his sudden death, which occurred in Castiglioncello on August 19, 1950. 
	
	\smallskip
	
	A member of the Academy of Italy from 1939 until the Academy's suppression, a National Fellow of the Lincei, a Pontifical Academician, 
	he belonged to numerous other Italian and foreign Academies and Associations.
	
	\smallskip
	
	Characteristic of Giovanni Giorgi was his interest in the most varied
	branches of science; his mind ranged from pure mathematics to the
	natural sciences (he willingly showed his rich collection of shells, classified with the most rigorous scientific method), to the most elevated questions of engineering. 
	
		\smallskip
	
	His breadth of vision is reflected in his publications (over 350)
	\footnote{
		A list of the publications of  Giovanni Giorgi can be found in the appendix of
		his book,  see Giorgi (1949)},
	among which one finds notes on mathematical physics alongside works of a strictly engineering nature, such as the project for the municipalization of the electricity plant of Avellino, or physical physiology papers, such as research on color vision.\\
	
	However, his major relevance lays on three topics:
	the system of units of measurement that bears his name (for which he  is universally known), the problem of absolute motion, 
	and  functional operational calculus.\\
	
	The Giorgi system of units is too well known to dwell on it here. It will suffice to remember that with it, the inconvenient electrostatic and electromagnetic systems are eliminated by establishing, in addition to the three mechanical units (meter, kilogram, second), a fourth unit of an electrical nature (as in the system electrostatic and electromagnetic systems in which the dielectric constant and the permeability of vacuum are set equal to one, respectively, but this leads to units that are too small or too large for common needs), taken from the practical system, for example, the ohm or the ampere. 
	For more details we refer the reader to the articles by the group of Arcangelo Rossi (e.g. Rossi et al. (2002), D'Agostino and Rossi (2004))
	and the most recent review by Frezza et al. (2015).
	His work in the second topic went almost unnoticed. We learned about it from the obituary
	of Giorgi by Dario Graffi (1951).
	The problem of absolute motion in the context of the fundamental laws of dynamics consists in finding the reference frame with respect to which the second law of dynamics is valid.
	Giorgi, appropriately interpreting the notion of force and the first law of dynamics, proves the second law to be valid with respect to any reference frame. This is consistent with the ideas of the general theory of relativity, which, however, was published, as everyone knows, in 1915, while Giorgi's memoir dates from 1912; indeed, he had been in possession of his results since 1903. The interest of Giorgi for the theory of relativity is proved by the brief correspondence with Einstein, see Frezza et al (2002). 
	
	\smallskip
	
	Even if Giorgi was interested mainly by electromagnetism and mathematical methods, his awareness of the 
	modern developments in physics is proven, for example, by his partecipation to the historical Como conference of 1927 in celebration of the centenary of Volta's death. In that conference Bohr introduced for the first time the principle of complementarity (see De Gregorio (2014)) and there was the partecipation of the fathers of quantum mechanics from Born to de Broglie, Pauli, Sommerfield, Planck and obviously the Fermi group (see Fig.1).
	
	\begin{figure}
		\centering
		\includegraphics[width=0.9\linewidth]{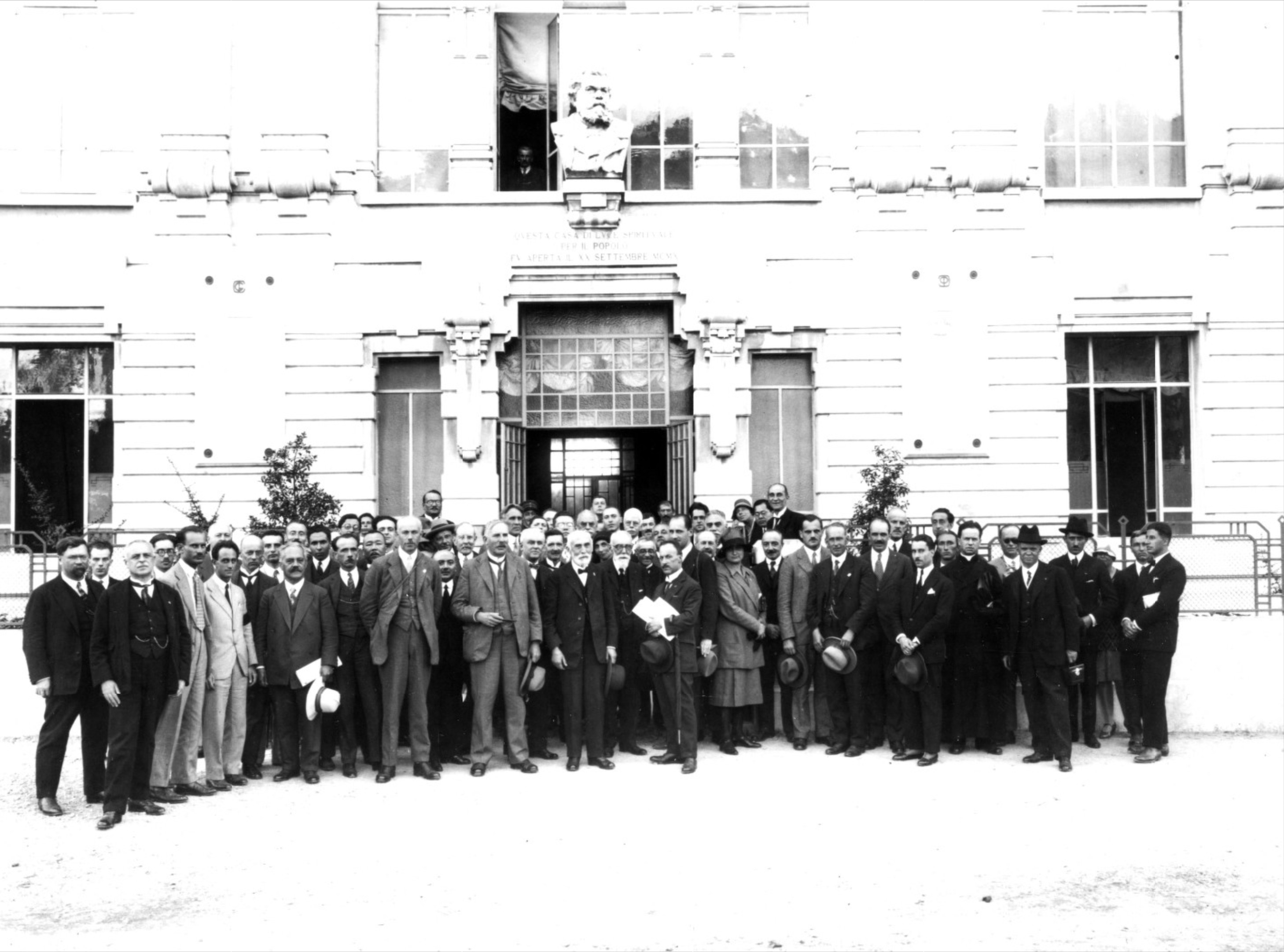}
		\caption{The photograph from the 1927 Como conference, held to commemorate the centenary of Volta's death. Professor Giorgi is clearly visible on the left, near Rasetti, Segrè, and Fermi. The source is the paper Goodstein (2001).}
	\end{figure}
	
	Giorgi was also a successful treatise writer; in addition to the volume mentioned above, his work on Rational Mechanics (Studium Urbis, Rome, 1946) and his lithographed lectures on Physics-Mathematics (Litografia Sampaolesi, Rome, 1928) were notable for their originality of views and clarity of exposition according to Graffi (1951).
	One must also mention his popular works (What is Electricity, Cremonese, Rome, 1928; Ether and Light, idem, 1939; The Fragmentation of the Atom, Ruffolo, Rome, 1946) and his excellent compendium of the History of Mathematics (S. E. L., Turin, 1948), a summary of a larger and still unpublished manuscript.
	
		\smallskip
		
	In this paper 
	we devote our attention to the third topic referring the reader to the next sections.
	On this respect we have mainly used the original papers of the author and the  obituary by Dario Graffi (1951) who, being a well-known professor of mathematical physics,
	mostly pointed out the mathematical aspects of the research activity of his colleague Giovanni Giorgi.\\
	Furthermore we added in an appendix the list of Giorgi's main papers 
	related to operational calculus.
	\\
	For the sake of completness, we should cite also the obituaries written by the celebrated mathematical physicist Signorini (1951a, 1951b),
	that however were not as helpful for this paper as the obituary by prof. Dario Graffi .\\
	Therefore, the novelty of this short paper is, in our view, twofold. First of all we point out the contribution of Giorgi to the mathematical studies
	about the operatorial methods to solve PDEs emerging in physical problems. In the review paper Giorgi (1934) shows a complete knowledge of the mathematica literature about 
this topic including contributions of well-known mathematicians such as Pincherle and Wiener. Secondly, we found that Giorgi had given a contribution also in the context of
 fractional differential calculus. Concluding, our main aim is to underline the original contribution of Giorgi in mathematical physics that is little known in the literature.

	\section{The contribution of Giorgi to the operational calculus}
	
	The great mathematical culture and original contributions of Giorgi in the context of the operational methods to solve mathematical problems 
	emerging in engineering and physics is proven by its high production in this field (see the Appendix). We consider, in particular, the paper Giorgi (1934), as a review paper that gives an historical view, an analysis of the main contributions given by the author and a collection of open problems. \\
	First of all, a relevant point is the close relationship between Giorgi's interest in that for this mathematical theory and the physical motivation lying in the analysis of electromagnetic problems. 
	As it is well known, functional operational calculus consists, in essence, in solving a partial differential equation or an integro-differential equation, in which the derivation with respect to a variable (the equation must be linear) can be treated as a numerical constant.
	The first, basic example, considered in Giorgi (1934) is the one-dimensional heat equation
	\begin{equation}
		\frac{\partial^2 U}{\partial x^2} = \alpha^2 \Delta U,
	\end{equation}
	where $\Delta = \partial/\partial t$ in the original notation of the Author. Then, the operational solution of the equation is given by (Giorgi, 1934)
	\begin{equation}\label{frac}
		U(x,t) = e^{-\alpha x\sqrt{\Delta}} U_1(t)+e^{\alpha x\sqrt{\Delta}} U_2(t).
	\end{equation}
	We observe that in the formal solution \eqref{frac}, there is a fractional derivative in time of order 1/2. This is the reason why Giovanni Giorgi devoted some studies to the meaning of derivatives of arbitrary real order. Most papers by Giorgi on operational methods are devoted to giving a rigorous mathematical meaning to this formal approach to solving linear PDEs in mathematical physics. \\
    Giorgi studied rigorously the mathematical literature on operational methods but his work was always motivated by concrete problems emerging in physical models.
    Another interesting example is given in Giorgi (1934) and was previously studied in detail in Giorgi (1928) regarding the signal trasmission in the telegraph. 
    The governing equation, in this case, is given by (in the notation of Giorgi (1934))
\begin{equation}
\frac{\partial^2 U}{\partial x^2} = a  \frac{\partial^2 U}{\partial t^2} + b \frac{\partial U}{\partial t} +c U,
\end{equation}
whose operational solution is represented as
\begin{equation}
	U(x,t) = e^{-\alpha x\sqrt{a\Delta^2+b\Delta+c}} U_1(t)+e^{\alpha x\sqrt{a\Delta^2+b\Delta+c}} U_2(t),
\end{equation}
 	where the coefficient $\Delta = \partial/\partial t$. Obviously this is a more complicated problem, see Giorgi (1928) for the complete treatment by operational method.\\
	In general, this procedure leads to finding the solution of
	the equation with expressions of the form
	$f\left(\frac{d}{dt}\right)\, V(t)$
	where $V (t)$ is a known function of $t$,
	$f \left(\frac{d}{dt}\right)$
	is a rational or transcendental analytic function of the operator considered
	as a numerical variable. 
	In order to give e precise meaning to the aforementioned expressions,  it is necessary, as they say, to evaluate the operator 
	$f \left(\frac{d}{dt}\right)$ applied to $V (t)$. Therefore, the basic problem is to give a rigorous meaning to this operation and the conditions for its application.\\
	
	The operational calculus had already been used, before Giorgi, by some
	authors, especially Heaviside, but only for particular problems, and the evaluation of  $f \left(\frac{d}{dt}\right)$
	was obtained with empirical methods,  as a mathematical trick to obtain exact solutions for differential equations emerging in electromagnetism. In some way, Heaviside was a pioneer in this field but he used correctly this method without a rigorous mathematical theory for the analysis of electromagnetic problems. Giorgi admired the work of Heaviside and he tried to give a rigorous mathematical theory for the symbolic methods developed by Heaviside.\\
	Recalling again Graffi's obituary, it is a merit of Giorgi
	to have legitimized the operational calculus by indicating a general
	rule for determining 
	$f \left(\frac{d}{dt}\right)\, V(t)$
	which thus results in a function $W(t)$ expressed
	by means of an appropriate generalization, due to Giorgi himself, 
	of the $Fourier$ integral (see Giorgi (1934)).
	Now, especially after the masterful works of Doetsch, the operational calculus is justified by the Laplace transform. That is, one
	first determines the transform of $W(t)$, 
	$\mathcal{L}[W(t)]$, which is then inverted,
	which can be done, for example, using the Fourier transform. 
	Therefore, the method of the
	Laplace transform and that of Giorgi differ very little in essence;
	in the first case one works on
	$\mathcal{L}[W(t)]$ itself,
	which is only inverted at the end of the calculations; in the other one, on 
	$W(t)$ itself, expressed by its generalized Fourier integral; that is, by means of the inverse transform of $\mathcal{L}[W(t)]$.\\
	
	There is another interesting historical point emerging from the paper of Giorgi (1934): the knowledge and appreciation of the parallel work developed in the same field by Norbert Wiener (1926) and the role of these operational methods in the quantum physics studies by Born and Dirac.
	This is again a proof of the deep knowledge of Giorgi of the recent developments in modern physics (first of all in relativity but also in quantum physics).

	\section{The contribution to the development of fractional calculus}
	The origin of fractional calculus (i.e. of the study of derivatives and integrals of arbitrary real order) goes back to the correspondence between
	L'H\^{o}pital and Leibniz concerning the meaning and relevance of derivatives of non-integer order. We refer to Ross (1977) for a first historical review. Also the first monograph devoted to the fractional calculus by Oldham and Spanier (1974) gave a good review of the main mathematicians that had contributed to that field. However, it is not trivial to re-construct the entire history of the Italian contribution to this field of research. We can refer, for example, to Mainardi and Pagnini (2011) about the role of the well-known mathematician Salvatore Pincherle in the development of fractional calculus and to Mainardi (2012) about the history of the relevant applications in the physical applications to viscoelasticity.\\
	
	\smallskip
	
	In our research we found that relevant contributions date back to other Italian mathematicians, including Giorgi, that in 1942 devoted a discussion paper to this topic. First of all, we would like to stress that in the paper Giorgi (1942), the Author considered, in particular, the papers written by Angela Maria Molinari about derivatives of arbitrary order. We think that some additional research about this figure should be done. It is well known that mathematical and physical studies were carried out mainly by males in that period. We found two not trivial and pioneering notes by A.M. Molinari (1916) on fractional derivatives published on Atti della Reale Accademia dei Lincei. \\
	As already seen, for example regarding equation \eqref{frac}, Giorgi approached the fractional calculus starting from the operational solutions of physical problems. The Author observed that arbitrary order derivatives in the general theory developed by Liouville and other authors were little known and used, due to an abstract approach, while an operational approach is strictly related to physical problems and can give a strong and unique theory.

	\section{Conclusions}
	In this note we remarked the role played by Giovanni Giorgi in mathematical physics by developing a rigorous mathematical theory for the operational methods used in physics by Heaviside and other authors for solving applied electromagnetic problems. We observe that, even today, these methods are object of research in applied mathematics, we refer for example to the work of Dattoli et al. (1997), Dattoli et al. (2004), Ricci and Tavkhelidze (2009) but there is a wide recent literature about this topic. Moreover, umbral methods are strictly related to the operational methods, we refer for example to the classical book by Roman and Rota (1978) and to the recent monograph Licciardi and Dattoli (2022).

		\subsection*{Acknowledgments}
	The authors are grateful to the librarians of the Department of Physics and Astronomy (DIFA) and of the Department of Mathematics of the University of Bologna for  providing   papers 
	by Giorgi published in the Proceedings of Academy of Lincei and of Academy of Italy.
	Furthermore they acknowledge the  availability to provide articles   on Giorgi
	published on  Physis, 
	RIVISTA INTERNAZIONALE DI STORIA DELLA SCIENZA, 
	(a journal edited for  the Societa' Italiana Storici della Fisica e 
	dell'Astronomia ,  SISFA)
	to 
	Dr Benedetta Campanile, University of Bari and to the members of the actual 
	Editorial Committee of Physis.
	
	The research activity of F. Mainardi
	has been carried out in the framework of the activities of the National Group of Mathematical Physics (GNFM, INdAM).

	\section{References}

	D'Agostino, S., Rossi, A. (Editors),
	Giovanni Giorgi nella realtà del suo tempo. 
	Physis, Vol 41 (2004), Special Issue
	
	\smallskip
	
	D'Agostino, S., Rossi, A.,
	Introduzione a Giovanni Giorgi nella realtà del suo tempo.
	Physis, Vol 41, pp. 309-317 (2004)
	
		\smallskip
		
Dattoli, G., Ottaviani, P. L., Torre, A., Vázquez, L. Evolution operator equations: integration with algebraic and finite difference methods. Applications to physical problems in classical and quantum mechanics and quantum field theory. La Rivista del Nuovo Cimento (1978-1999), 20(2), 3-133, (1997)

	\smallskip
	
	Dattoli, G., Ricci, P. E., Khomasuridze, I., Operational methods, special polynomial and functions and solution of partial differential equations. Integral Transforms and Special Functions, 15(4), 309-321 (2004)
	
	\smallskip
	
	De Gregorio, A., Bohr's way to defining complementarity. Studies in History and Philosophy of Science Part B: Studies in History and Philosophy of Modern Physics, 45, 72-82, (2014)
	
	\smallskip
	
	Frezza, F.,  et al.,
	The Life and Work of Giovanni Giorgi: The Rationalization of the International System of Units.
	IEEE Antennas and Propagation Magazine,  Vol. 57, No 6, December (2015),
	pp. 152--165.
	DOI: 10.1109/MAP.2015.248675
	
	\smallskip
	
	Giorgi, G., Sugli integrali dell'equazione di propagazione in una dimensione.
	Rendiconti del Circolo Matematico di Palermo, Vol. 52, pp. 265-312, (1928)
	
	\smallskip
	
	Giorgi, G., Metodi moderni di calcolo operatorio funzionale. 
	(Conferenza tenuta al Seminario Matematico e Fisico di Milano, aprile 1934). Rendiconti del Seminario detto, Vol. VIII, pp. 189-214, (1934).
	
	\smallskip
	
	Giorgi, G., Formole per la derivazione a indice generalizzato. 
	Rendiconti R. Accad. d'Italia, seric 7,vol. III, pp. 693-701 (1942).
	
	\smallskip
	
	Giorgi, G., Verso l'Elettrotecnica Moderna: Richiami e Contributi. Tamburini, Milano, (1949).
	
	\smallskip
	
	Goodstein, J. R., A conversation with Franco Rasetti. Physics in Perspective, 3(3), 271-313, (2001).
	
	\smallskip
	
	Graffi, D. ,
	Necrologio di Giovanni Giorgi.
	Bollettino dell’Unione Matematica Italiana, 
	Serie 3, Vol. 6, n.2, pp. 171--188, (1951).
	Zanichelli, Bologna.
	
	\smallskip
	
	Licciardi, S., Dattoli, G., Guide to the Umbral Calculus: a Different Mathematical Language, World Scientific, Singapore,  (2022).
	
	\smallskip
	
	Mainardi, F., An historical perspective on fractional calculus in linear viscoelasticity: short survey. Fractional Calculus and Applied Analysis, 15(4), 712-717, (2012).
	
	\smallskip
	
	Mainardi, F., Pagnini, G., The role of Salvatore Pincherle in the development of fractional calculus. In
	S. Coen (editor), Mathematicians in Bologna 1861–1960 (pp. 373-381). Basel: Springer Basel, (2011).
	
	\smallskip
	
	Molinari, A.M., Derivazione ad indice qualunque, Rendiconti Lincei, Serie V, vol. 25, pag. 230-233 and pag. 268-273, (1916-B).
	
	\smallskip
	
	Oldham, K., Spanier, J., The Fractional Calculus, Theory and Applications of Differentiation and Integration to Arbitrary Order (Vol. 111). Elsevier, (1974).
	
	\smallskip
	
	Ricci, P. E., Tavkhelidze, I., An introduction to operational techniques and special polynomials. Journal of Mathematical Sciences, 157(1), 161-189, (2009).
	
	\smallskip
	
	Roman, S. M. and Rota, G. C., The Umbral calculus. Advances in Mathematics, 27(2), 95-188, (1978).
	
	\smallskip
	
	Ross, B., The development of fractional calculus 1695–1900. Historia mathematica, 4(1), 75-89, (1977).
	
	\smallskip
	
	Rossi, P.,
	Cento anni fa il primo concorso italiano di fisica teorica.
	Giornale di Fisica (SIF), Vol 67. No. 1, pp 53--60.
	Societ\`{a}Italiana di Fisica (SIF), Bologna.
	DOI 10.1393/gdf/i2026-10636-4, (2026)
	
	\smallskip
	
	Rossi, A., D'Agostino, S., Morando A.P.,
	Giovanni Giorgi e la tradizione dell'Elettrotecnica Italiana
	Societa' Italiana Storici della Fisica e dell'Astronomia (SISFA)
	In: Atti del XXII Congresso Nazionale di Storia della Fisica e dell’Astronomia, Genova e Chiavari, 6-8 giugno 2002, p. 314-328, (2002).
	\\
	https://www.sisfa.org/pubblicazioni/atti-del-xxii-convegno-sisfa-genova-chiavari-200
	
	\smallskip
	
	Signorini, A., In Memoria di Giovanni Giorgi, Il Nuovo Cimento, Volume 8, pages 225–228, (1951a).
	
	\smallskip
	Signorini, A., 
	Commemorazione del Socio Giovanni Giorgi,
	Rendiconti Lincei,  Vol. 11 , No 6 (1951), pages 4126--423, (1951b).

	\smallskip
	Wiener, N., The Operational Calculus, Mathematische Annalen, 
	95, pp. 537-583, (1926).

	\newpage
	\section{Appendix: List of the Giorgi's papers devoted to operational methods and fractional derivatives}
	
	G.Giorgi, II metodo simbolico nello studio delle correnti variabili.
	Comunicazione al Congresso Elettrotecnico di Napoli, ottobre 1903. 
	Atti  A.E.I., vol. VIII, 1904, pp. 65-43.
	
	\smallskip 
	
	G. Giorgi, Sul calcolo delle soluzioni funzionali originate dai problemi di elettrodinamica. 
	Comunicazione dell'A.E.I. Roma, 26 giugno 1905. Atti A.E.I., vol. IX, 1905, pp. 651-699.
	
	\smallskip
	
	G. Giorgi, Sulla teoria delle equazioni integrali e delle loro generalizzate. 
	R. Accademia Lincei. Rendiconti, vol. XXI (1° sem. 1912), pagine 748-754.
	
	\smallskip
	
	G. Giorgi, Sui problemi dell'elasticit\`a ereditaria. 
	R. Accad. Lincei, Rendiconti, vol. XXI (2° sem. 1912), pp. 412-418.
	
	\smallskip
	
	G. Giorgi, Sugli operatori funzionali ereditari. 
	R. Accademia Lincei, Rendiconti, vol. XXI (2° sem. 1912), pp. 683-687.
	
	\smallskip
	
	G. Giorgi, Sugli integrali dell'equazione di propagazione in una dimensione.
	Rendiconti del Circolo Matematico di Palermo, Vol. 52, 1928, pp. 265-312.
	
	\smallskip
	
	G. Giorgi, Sulla propagazione delle onde nei mezzi con assorbimento selettivo. R. Accademia Lincei, Rendiconti, serie 6,vol. IX (gennaio 1929), pp. 8-11.
	
	\smallskip
	
	G. Giorgi, Metodi per la calcolazione dei fenomeni transitori nel regime variabile delle correnti. 
	Dati e Memorie sulle Radiocomunicazioni, pubblicati dal Consiglio Nazionale delle Ricerche, vol. IV, 1932, pp. 581-589.
	
	\smallskip
	
	G. Giorgi, Metodi moderni di calcolo operatorio funzionale. 
	(Conferenza tenuta al Seminario Matematico e Fisico di Milano, aprile 1934). Rendiconti del Seminario detto, Vol. VIII 1934, pp. 189-214.
	
	\smallskip
	
	G. Giorgi, Metodi di calcolo per la propagazione dei segnali nelle linee. 
	Rendiconti del Congresso 1936 dell'A.E.I. in Roma. Atti del Congresso, 1937.
	
	\smallskip
	
	G. Giorgi, Questioni sul calcolo operatorio funzionale, 
	in Conferenze di Matematica e Fisica, della R. Università di Torino, vol. 1934-35-36, pp. 211-222. (Conferenza tenuta il 27 aprile 1936).
	
	\smallskip
	
	G. Giorgi, Metodo pel calcolo degli effetti di distorsione nelle linee telegrafiche e telefoniche. 
	R. Accademia Lincei. Rendiconti, serie 6, vol. XXV (febbraio 1937), pp. 155-156.
	
	\smallskip
	
	G. Giorgi, Calcolo dei fenomeni transitori nei circuiti a corrente alternata. 
	R. Accademia Lincei, Rendiconti, vol. XXIX, serie 6 (maggio 1939), pp. 541-542.
	
	\smallskip
	
	G. Giorgi, A proposito di una recente pubblicazione sul calcolo operatorio. Elettrotecnica, vol. XXVI I, 940, pp. 528-529.
	
	\smallskip
	
	G. Giorgi, Sulla funzione impulsiva nel calcolo operazionale.
	Lettera alla Redazione, L’Elettrotecnica 10 maggio 1941, vol. XXVII, n. 9, p, 229.
	
	\smallskip
	
	G. Giorgi, Formole per la derivazione a indice generalizzato. 
	Rendiconti R. Accad. d'Italia, seric 7,vol. III (novembre 1942), pp. 693-701.
	
	\smallskip
	
	G. Giorgi, Progressi e applicazioni del calcolo operatorio funzionale. 
	Atti del Convegno Volta della R. Accad. d'ltalia, 1943, pp. 281-290

\end{document}